\documentstyle{amsppt}
\magnification=1200
\nologo

\def\ms{\medskip}

\def\T{\Theta}

\def\oonL{\mathop{{\hbox{$L^{m,2}$}\kern -20pt\raise7pt
\hbox{$\circ$}}}}

\def\O{\Omega}

\def\emp{\emptyset}
\def\bigtimes{\mathop{{\lower2.95pt\hbox{$\wedge$}\kern-6.666pt\raise2.95pt
\hbox{$\vee$}}}\limits}
\def\crd{\cr\noalign{\vskip4pt}}


\def\eq{\eqalign}

\def\e{\epsilon}

\def\l{\lambda}
\def\R{${\text{\bf R}}^N$}
\def\no{\noindent}
\def\bs{\bigskip}
\def\G{\Gamma}
\def\g{\gamma}
\def\p{^\prime}

\def\b{\beta}
\def\a{\alpha}
\def\n{\nabla}

\def\da{d_{\partial \O}(x)}
\def\d{\da}
\def\pn{\par\noindent}


\hsize=5in
\vsize=7.6in
\baselineskip=16pt

\def\Wmps0O{W^{m,p}_0(\O,\d^s)}
\def\Wmpt0O{W^{m,p}_0(\O,\d^t)}

\NoRunningHeads
\topmatter
\title
Polynomial Capacities, Poincar\'e type Inequalities
and Spectral Synthesis in Sobolev Space
\endtitle
\footnote""{AMS classification numbers: 26D15, 46E35 and 35J30.}
\author
Andreas Wannebo
\endauthor
\address
Department of Mathematics
\newline The Royal Institute of Technology
\newline
Stockholm
\newline
Sweden
\endaddress
\abstract
There are polynomial capacities defined by Maz'ya, that solves the 
problem of finding (approximately) the constant in a Poincar\'e type 
inequality in a cube. 
A similar inequality has been studied by Hedberg and we define
a new set of polynomial capacities that gives the correct answers to the 
corresponding problem of finding the constant in the generalized
setting of the Hedberg case. 
\par
All this is a prerequisite for a planned publication on Hardy inequalities. 
\par
Applications are given here that also show the importance of these ideas. 
\par
This concerns the Sobolev space Spectral Synthesis problem solved by
Hedberg and also treated by Netrusov.
An equivalent formulation in terms of polynomial capacities is given
proved by new and important technique. This technique is also applied to give a
concrete new result similar to Spectral Synthesis in Sobolov space,
but for the nonnegative cone in Sobolev space of order 2 instead.
\endabstract
\endtopmatter

\document

\heading
Introduction I.
\endheading

The Sobolev spaces with integer order of derivatives are often
seen as a just a special kind of Bessel potential spaces, since their
norms generally are equivalent in the analysis sense of the word.

We will see here the importance of the point of view of
Sobolev space norm structure as a sum of seminorms (or norm equivalent to such a
norm). We want to emphasize this point.

Having for instance a set of Sobolev functions, then the relationships
between these seminorms in terms of inequalities are of key importance.

This we will show in the last part of this paper and in a coming
paper on Hardy inequalities. We refer also to [WAN] were this
connection where exploited.
\bs

The material here as well as the forthcoming paper on Hardy
inequalities has been widely circulated in some different versions.
The first not that very circulated (typed) version dates back to the
early 80:ies.
\bs

For technical reasons these relationsships will be studied with a unit
cube as the domain. The cubes give proper versiality e.g. by cubic
decompositions like the Whithney cube decompositon and other cube
decompositions.

As this is an introduction only we do not always give precise definitions,
results etc. We instead refer to the main body of the text.
\bs

Now a rough description of the Poincar\'e (type) inequalities in
question.
\bs

We simplify notation and let temporarily $|||u|||$ denote
the $L^p$-norm of the function $u$ from a unit cube $Q$ in \R\ to 
$\text{\bf R}$.
We put $D^{\a}u$ to be the weak partial derivative of $u$ with a 
multiindex for the derivative.

Furthermore we put 

$$
|||\n^ku|||=\sum_{|\a|=k}\ |||D^{\a}u|||
$$
and
$$
|||u|||_W=\sum_{k=0}^m\ |||\n^ku|||
$$
as the corresponding Sobolev space norm with parameters $m,p$. 

This formulation emphasizes the sum of seminorms structure. Although
the definition can be made differently, this is of no real effect here.
The Sobolev space $W$ belonging to this norm
is the set of all functions (or more properly distributions) that makes the
Sobolev space norm exist.
\bs

Now the Poincar\'e type inequality in a cube can be approached as
follows. 
\bs

Question A: For all $u$ in $\Cal A$, an arbitrary subset of $W$, study the following 
inequality,

$$
\sum_{i=0}^m b_i|||\n^iu|||\le \sum_{i=0}^m a_i|||\n^iu|||.
\leqno(1.1)
$$

What can then be said about the possible sets of constants $\{a_i\},\{b_i\}$?

This a quite general formulation.
\bs

The way the Question A is approached here is much simpler though.

First these constants are made dependent of one parameter only. 
This can be written in such a way that one of the constants is
selected as the parameter, (the choice clearly makes a difference).
Hence the remaining constants are functions of the one selected.

Now the Question A reduces to find this selected constant as a 
function of $\Cal A$.
\bs

There are simple intrinsic relations - interpolation inequalities -
to the effect that e.g. the intermediate gradient seminorm is less than or
equal to the sum of a lower and a higher gradient seminorm multiplied
with a constant. 

This simplifies the possibilities, since our analysis of
the Question A is only up to a constant -- a constant which however can be 
approximately calculated.
\bs

So far we know of only two cases of the Question A being treated,
modulo interpolation i.e.
\bs

The first case was given by Maz'ya in [MAZ]. He let $\Cal A$
unnecessarily be a vector space, but all works out well with a
general subset. All that is needed is to change definitions
somewhat and to follow this up in the proof.

This smallish distinction happens to be of great importance as will be seen
here.
\bs

The second case is quite similar and given by Hedberg see [HED],
He however did not treat the question the way Maz'ja did. He wanted
the result as a lemma for his proof of Spectral synthesis in Sobolev
space see [HED] or better [A-H].
Hence he lost the generality of the question and also
he did only got a weaker result with ordinary non-polynomial capacity
for his special case. -- He wanted a lemma.
\bs

The Maz'ya case

$$
\sum_{i=0}^k |||\n^iu|||\le C\ (|||\n^{k+1}u|||+|||\n^mu|||).
\leqno(1.2)
$$

Here Maz'ya proves that the constant $C$ is equivalent to the quantity $\G$
to power $-1/p$, (with proper parameters), as function of $\Cal A$,
(as above).
In order for the equivalence to hold it is enough that $\G$ is smaller
than some constant or that $C$ is larger than a constant.

Maz'ya has also treated the situation with small values of $C$ etc in [MAZ],
but is not of interest in the applications given here.

This situation also has a bearing on our way to do Hardy inequalities but
the instance with large $C$ is much more important and is what will be
used there.
\bs

The Hedberg case can be formulated,

$$
\sum_{i=0}^k |||\n^iu|||\le C_0\ |||\n^{k+1}u|||+C|||\n^mu|||),
\leqno(1.3)
$$
with $k+1$ less than or equal to $m$, and where $C_0$ is a fixed
constant and where the problem is to find the constant
$C$ as a function of $\Cal A$.
\bs

Now it happens that solving this two instances of the Question A solely 
depends on a certain inequality. Here named the weak Poincar\'e inequality.

It states that given $k$ and $p$ there is a constant $A$ such that
for every $u$ with RHS finite there exists a polynomial of degree
less than or equal to $k-1$ such that

$$
\sum_{i=0}^{k-1}\ |||\n^i(u-P)|||\le A|||\n^ku|||.
\leqno(1.4)
$$

If $p$ is finite and bigger than one this can be done as
a unique projection of $u$ to the polynomials of degree less than $k$
by taking an infimum. This is possible because of local convexity.
\bs

Maz'ya in his definition makes some extra work on his definition for 
pedagogical(?) reasons. Say $\Cal A$ are the Sobolev functions
(with the relevant parameters), which are the set formed by the
restriction of to the cube of the closure of
$C^{\infty}_0(K^c)$, with $K$ compact subset of the closure of the
cube. Now $\G$ as function of $\Cal A$ can in fact be seen as a
function of $K$. That it is seen as function of sets of points in \R\ or 
rather in the cube. The reason for chosing the special appearence of
$\G$ with that inverted $p$-th root will now be made clear.

If $k+1$ equals 1 then 
$\G$ is a quantity very close to the definition of a capacity for
Sobolev space in a so called condenser formulation and with 
parameters $m$, $p$.
 However it is possible to prove that $\G$ here is in fact
equivalent to this capacity and can serve as an alternative definition.
 For higher values of $k+1$ the indicator functions $1_K$ is replaced by
polynomials of higher degree times the indicator functions $1_K$.
To get $\G$ involves however an infimum over these possible polynomial
as well, roughly speaking. This has the consequence that subadditivity
fails without anything to replace it with for $\G$.
\bs

This is the reason for naming these quantities polynomial capacities.
\bs

We treat the generalization of the Hedberg case solving the question
completely and but with methods that are completely different from his.

Hedberg case we can say that not only do we
solve the question completely and with methods radically different, but
also the formulation of the solution is also radically different.
\bs

Now back to the formulations and proofs.
It clear to any kean observer that in
proving general statements like this there is something tautological
to the proof process. In this case it is so too, but modulo the weak
Poincar\'e inequality.

Now there is an important distinction in mathematics between bad
tautologies and good tautologies. A bad one takes you from something
you do not know to something you do not know. A good tautology
takes you from something you do not know to something you can know
about, maybe know a lot about.

The situation at hand here is that of good tautologies and the results 
are very useful.
\bs

We claim that the concepts of polynomial capacities is
fundamental to the theory of Sobolev spaces.
\bs

Summarizing:
\bs

I. The Question A, at least the cases treated here, is of a fundamental
nature, i.e. the polynomial capacities turn up in many places where
they primarily were not really asked for.
\ms

II. The polynomial capacities have interesting properties, 
e.g. in the
discussion above on the polynomial nature of these creatures, when
looked upon as functions of point sets, then in fact the polynomial
capacity measures quantatively how much the point set deviates from 
belonging to any algebraic surface within a certain class (this is
somewhat simplified though).
\ms

III. The polynomial capacities can often in practise be estimated in 
good ways. Hence there can be concrete answers to problems where
they turn up as tools.
\vfill\eject

Now we turn back to a discussion of the case in the general setting.

The technicalities in the
definition of the polynomial capacity $\T$ that takes the place of
the Maz'ya polynomial capacity $\G$ are more difficult. In spite of
this many tries for good estimates have turned out rather
well. This after harder work i.e. One can say that as long as the
polynomial capacity $\T$ has the same estimate as the polynomial
capacity $\G$, the second Poincar\'e inequality in fact provides a
stronger result. This is of course obvious. 

What is not obvious that we do not know instances where these
polynomial capacities really are different. Hence we here have an open
question.
\bs

The situation that the set $\Cal A$ 
is not a vector space turns up in
the treatment of $\T$ for $\Cal A$ a nonnegative cone. The 
implication is that nonnegative functions sometimes have smaller (=
better) constant
$C$ in the inequality than some related set of general functions.
This key result will used be for the last part of this paper as well as
be very important in the coming study of Hardy inequalities.
\bs

Now we have discussed the study of the polynomial capacity part. There
is only to add that in our treatment the theorems are furnished with a
more complicated dress of different seminorms. This can be ignored at
first reading.
\bs

Next we discuss the applications of the polyonomial capacities made here.
\bs

This makes up for an achievement that should rather been seen as a new
and intuitively simple technique made possible by the apparence of 
these polynomial capacities and a delicate use of them.
\bs

Say we have a Sobolev space $W'$, with parameters $m,p$, and defined
in all \R. Say furthermore that we have two closed subspaces say,
$F_1$ and $F_2$ of $W'$. 
\bs

A question here can be: Is e.g. $F_1$ a subset of $F_2$.
\bs

An answer to such a question in a more specific situation, yet still in 
extreme generality, is really what makes up the proof of the important
theorem called Spectral synthesis in Sobolev space proved by 
Hedberg [HED] or better in [A-H] and later by Netrusov see [A-H]. 
The latter proved this in a
setting of more general function spaces.
\bs

The reader is refered to this the last part of this paper for the
details since these are to lenghty to suit the introduction.
\bs

Now we can in a somewhat abstract way discuss at least the ideas
behind this new technique. 

Say that $F_2$ is a subset of $F_1$. Say the subspaces have the
property that there is closed set $K$ such that the subspaces are 
the exactly the same away from $K$, i.e. for any compact in the complement
of $K$ the subspaces coincides. Hence
what may make $F_1$ not be a subset of $F_2$ is concentrated to 
any neighborhood of $K$.
\bs

Proceed by chosing any $u$ in $F_1$.
The idea now is to make a $\e$-sized cubic decomposition of \R.
Form the set of cubes that intersects $K$. Label it $\{Q'\}$.
Then to start with take a cube $Q'$ in $\{Q'\}$.

Next take the restriction of $u$ to $Q'$ and construct another function
$u'$ only defined in $Q'$, which has the right boundary values to
fit as a replacement of $u$ in $Q'$, i.e. $u_1$ still belongs to
Sobolev space $W'$. Outside the union of $\{Q'\}$ it of course
belongs to $F_2$. However there is another property that $u'$ shall
have namely in the cube $1/2Q'$ (same centre and orientation but half
the side length) $u'$ shall coincide with a function in $F_2$.
\bs

Here is where the polynomial capacities come into the picture.
There is a lemma that that proves that this replacement possible if:
The polynomial capacities for the two subspaces are equivalent. 
(Actually this is stated more detailed.)
\bs

Next we iterate this enough to have completely changed the $u$ in $F_1$
to a $F_2$ function. Then there is the question of the cost in norm
for all these changes. In the situations studied here they go to zero
with $\e$. Hence that $F_1$ is a subspace of $F_2$ then follows from 
that these subspaces are closed.
\bs

{\bf Formulation of Spectral Synthesis in Sobolev space A:} It is the
result that, given a Sobolev space $W'$, then for very closed subset
$K$ of \R\ the subspace defined as the subset of all functions that have
all traces equal to zero for every partial derivative of order less 
than $m$ on this closed set $K$ coicides with to the subspace defined as the
the closure in the $W'$ norm of the infintely differentiable functions
with compact support outside $K$.
\bs

{\bf New Result B:} Spectral Synthesis in Sobolev space is
equivalent to the fact that certain related polynomial capacities are
equivalent.
\bs

{\bf New Result C:} Take $W'$ with $m$ equal to 2. Take its nonnegative
cone. For any closed set $K$ in \R\ the set of functions with partial 
derivative of order zero, i.e. no differentation, with zero trace on
$K$ is equal the closure in the $W'$
norm of the infinitely differentiable functions with compact support 
outside $K$.
\bs

This new result gives the first "not terribly complicated''
proof in this buissness with a $m$ greater than 1. (The $m$ equal to 1 
case is simple is this context).
\bs

After this descriptive and popularized introduction we turn to the
main body with all the gory details.
\vfill\eject

\heading
Section 2. On Polynomial Capacities.
\endheading
\bs

The polynomial capacities were invented by Maz'ya, see [MAZ]. He studied
a special kind of Poincar\'e inequality in a cube $Q$. The cube can
be taken as unit size since the general case is got from a dilation
of this case. The advantedge is of course that the formulas get rid of
a lot of factors that are only introduced by the scaling and not
really have interest by their own.
\bs

Other domains than cubes can be studied but if they are a bit bad
then nothing is known about the polynomial capacities.
\bs

The polynomial capacities in the
Maz'ya sense were reinvented by Bagby in a special dual setting. (This
pointed out by me - Use duality and the Hahn-Banach theorem.)

He gives a long and technically interesting paper however since he
makes use of arbitrary open sets instead of cubes for the polynomial
capacities. There seem to be no way verify his conditions in the main
theorem. Netrusov has made a follow up of this paper.
\bs

Now we state the Poincar\'e inequality studied by Maz'ya,
see [MAZ; 10.3].

The presentation here is a bit different, but there is in fact no
major change. The differences are that the polynomial capacities
are defined as function of a ``general'' set of functions instead of a
vector space and that more general set of semi norms (quasi-semi
norms) are possible.
\bs

First some preliminaries.
\bs

\no{\smc 2.0. Definition.} Let $\O$ be an open set in \R\ and $u$ a 
distribution from $\O$ to $\text{\bf R}$ then if the $D^{\a}u$ is 
in $L^p(\O)$ for all $\a$ of order less than or equal to $m$ then
$u$ is said to belong to the Sobolev space $W^{m,p}(\O)$, a Banach space,
and the definition of the corresponding norm is

$$
||u||_{W^{m,p}(\O)}=\sum_{k=0}^m||\n^ku||_{L^p(\O)},
$$
where
$$
||\n^ku||_{L^p(\O)}=\sum_{|\a|=k} ||D^{\a}u||_{L^p(\O)}.
$$
\bs

\no{\smc 2.1. Notation.} Let ${\Cal P}_k$ denote the polynomials of degree
less or equal to $k$ in \R.
\bs

\no
{\smc 2.2. Definition.} Let $P$ be a polynomial, let $Q$ be a fixed cube in \R,
let ${\Cal A}\subset W^{m,p}(Q)$ and let $1\le p$. Denote

$$
U_{P,{\Cal A}}=\{u\in W^{m,p}_0(2Q):u-P|_Q\in {\Cal A}\}.
$$
Let
$$
\G_{m,k,p}({\Cal A})=\inf_{P\in {\Cal P}_k}\ \inf_{u\in U_{P,{\Cal A}}} 
{||\n^mu||_{L^p}^p\over ||P||_{L^p(Q)}^p}.
$$
\bs

\no{\smc 2.3. Definition.} Define $C^{m,\l}(\O)$, H\"older space, 
as the subspace of $C^m(\O)$, which consists of the functions $u$ such that
the norm

$$
||u||_{C^{m,\l}(\O)}=\sup_{0\le i\le m}\ ||\n^iu||_{L^{\infty}}
+
\sup _{x,y\in\O,|\a|=m}
\ \{{|D^{\a}u(x)-D^{\a}u(y)|\over|x-y|^{\l}}\}
$$
is finite.
\bs

Now the Maz'ya theorem for his Poincar\'e inequality but in our presentation.
\bs 

\proclaim{2.4. Theorem}
Let $Q$ be a fixed cube in \R.
Let $T$ be one of the spaces below:
\roster
\item"(i)"
$L^q(Q)$ with $0<q\le {Np/(N-mp)}$ when $N> mp$ or $p=1$ and $N=m$;
\item"(ii)"
$L^q(Q)$ with $0<q<\infty$ when $N=mp$;
\item"(iii)"
$L^q(Q)$ with $0<q\le\infty$ when $N<mp$;
\item "(iv)" 
$C^{h,\l}(Q)$ for $(m-h)p>N>(m-h-1)p$ with $0<\l\le m-h-{N/p}$;
\item"(v)" 
$C^{h,\l}(Q)$ for $N=(m-h-1)p$ with $0<\l<1$.
\endroster
\par

Let ${\Cal A}\subset W^{m,p}(Q)$, let $1\le p$ and let $0<p_0$. 
Then the inequality
$$
||u||_T\le C\left(||\n^{k+1}u||_{L^{p_0}(Q)}
+||\n^mu||_{L^p(Q)}\right)\leqno(2.0)
$$
is valid for some constant $C$ iff $\G_{m,k,p}({\Cal A})>0$.

Furthermore the constant $C$ is equivalent to 
$(\G_{m,k,p}{\Cal A})^{-1\over p}$ for $C$ large or 
$\G_{m,k,p}({\Cal A})$ small.

\endproclaim
\bs

\demo{Proof}
First we note that the normalization made by Maz'ya's version is unnecessary.
Hence $\Cal A$ can be a general subset instead of a vector subspace.

Next we have by the Sobolev imbedding theorem and the Maz'ya result
i.e. his original result that
$$
\eq{
||u||_T\le &A\left(||u||_{L^p(Q)}+||\n^mu||_{L^p(Q)}\right)
\crd
\le 
&A(\G_{m,k,p}({\Cal A})^{-{1\over p}}
\sum_{j=k+1}^m||\n^ju||_{L^p(Q)}
+||\n^mu||_{L^p(Q)}).
\cr}\leqno(2.1)
$$

First, if $p_0\ge p$, we have by interpolation, see
[RAD; 4.13], that
$$
||u||_T
\le 
A\left(\left(\G_{m,k,p}({\Cal A})\right)^{-{1\over p}}
\left(||\n^{k+1}u||_{L^p(Q)}+||\n^mu||_{L^p(Q)}\right)
+||\n^mu||_{L^p(Q)}\right)
$$
and the result follows by the H\"older inequality.

Now let $p_0<p$.
There is a theorem by Sobolev, see [MAZ, 1.1.15], with as
a special case that for $|\a|=k+1$ and $p_0<p$ we have that
$$
\sum_{j=0}^{m-k-1}||\n^jD^{\a}u||_{L^p(Q)}
\sim||D^{\a}u||_{L^{p_0}(Q)}+||\n^{m-k-1}D^{\a}u||_{L^p(Q)}.
$$

Hence

$$
\sum_{j=k+1}^m||\n^ju||_{L^p(Q)}
\sim ||\n^{k+1}u||_{L^{p_0}(Q)}+||\n^mu||_{L^p(Q)}.
$$

This combined with (2.1) gives

$$
||u||_T
\le 
A\left(\left(\G_{m,k,p}({\Cal A})\right)^{-{1\over p}}
\left(||\n^{k+1}u||_{L^{p_0}(Q)}+||\n^mu||_{L^p(Q)}\right)+
A||\n^mu||_{L^p(Q)}\right).
$$

End of proof.
\bs

Next we turn to the Poincar\'e inequality devised by Hedberg,
which we will do very differently and in way completely.

The idea here is again to define an appropriate polynomial capacity 
to describe the second constant $C$, ($C_0$ is fixed). 

Hedberg, see [HED] or [A-H], used a Bessel capacity and got a sufficient condition only,
as well as for somewhat special situation.
\bs

We define a well-known projection operator from Sobolev space
functions to polynomials.  
\bs

\no{\smc 2.5. Definition.} Let $Q$ be a fixed cube in \R\ and let $u\in
W^{r+1,p}(Q)$, with $1\le p$.
Let $\Pi_{r,r,p}u$ denote a $P$ giving minimum (if $p>1$, at least, it exists
uniquely) in 
$$
\inf_{P\in {\Cal P}_r}\sum_{i=0}^r||\n^i(u-P)||_{L^p(Q)}.
$$
\bs

The fundamental property of the projection operator is the
following well-known inequality, see [MAZ; 1.1.11], not
formulated in the most general way. 
\bs

We denote this inequality {\bf the weak Poincar\'e inequality.} 
\bs

Warning we will constantly refer to it by this name!

(The naming of inequalities is a notorious matter.)
\bs

\proclaim {2.6. Theorem}  Let $\Omega$  be open, bounded in \R,
connected and have the cone property. Let $u\in W^{r+1,p}(\Omega)$ and
let $1\le p$, then

$$
\sum_{k=0}^r||\nabla ^k (u-\Pi_{r,r,p}u)||_{L^p(\O)}
\le A_\Omega||\n^{r+1}u||_{L^p(\Omega)},
$$
where $A_{\O}$ depends on $\O$, $r$ and $p$.
\endproclaim
\bs

Another source for this kind of inequalities is, [MEY], by N.G. Meyers.
\bs

Next we again map the polynomial $\Pi_{r,r,p}u$ to several others
that we will use frequently.
\bs

\no{\smc 2.7. Notation.} Let $Q$ be a fixed cube in \R\ and let $u\in
W^{r+1,p}(Q)$, with $1\le p$.
Let $0\le r'\le r$, let 
$\Pi_{r,r',p}u$ denote the $r'$-degree part of the polynomial
$\Pi_{r,r,p}u$ and
let $\Pi^{\text{\rm co}}_{r,r',p}u$ denote the complementing
polynomial. Furthermore denote

$$
\Pi^{\p}_{r,r',p}u=\Pi_{r,r',p}u-\Pi_{r,0,p}u.
$$
\bs

Now this projection operator is used to define a new kind of
polynomial capacity.

Remember that the parameter $\a$ below has a curious role to play. 
It only
is there to ``neatify'' the process later and do not really play the 
role of an ordinary parameter.
When thinking about this matter $\a$ can be ignored and hence the
situation is somewhat less messy.
\bs

\no{\smc 2.8. Definition.} Let $Q$ be a fixed cube in \R. Let $p\ge 1$,
let $\a>0$ and let ${\Cal A}\subset W^{m,p}(Q)$. Define

$$
\eq{
T^{\a}_{m,k,p,{\Cal A}}
=
\big\{
u\in {\Cal A}:
||\Pi_{m-1,k,p}u||_{L^p(Q)}
\ge
&\a||\Pi^{\text{\rm co}}_{m-1,k,p}u||_{L^p(Q)}
\crd
&\text{and}
\crd
||u-\Pi_{m-1,m-1,p}u||_{L^p(Q)}\le&{1\over 2}||\Pi_{m-1,m-1,p}u||_{L^p(Q)}
\big\}
\cr}
$$

and

$$
\T^{\a}_{m,k,p}({\Cal A})
=\min \left[\inf_{u\in T^{\a}_{m,k,p,{\Cal A}}}
{||\n^mu||^p_{L^p(Q)}\over ||\Pi_{m-1,k,p} u||^p_{L^p(Q)}},1\right].
$$
\bs

The capacity $\T^{\a}_{m,k,p}$ is taylored for the following
Poincar\'e inequality \'a la Hedberg.
\bs

\proclaim{2.9. Theorem}
Let $Q$ be a fixed cube in \R, 
let ${\Cal A}\subset W^{m,p}(Q)$, let $1\le p$ and let $0<p_0$. Let $T$ be 

\roster
\item"(i)"
$L^q(Q)$ with $0<q\le {Np/(N-mp)}$ when $N>mp$ or $p=1$
and $N=m$;
\item"(ii)"
$L^q(Q)$ with $0<q<\infty$ when $N=mp$ and $p>1$;
\item"(iii)"
$L^q(Q)$ with $0<q\le\infty$ when $N<mp$;
\item "(iv)" $C^{h,\l}(Q)$ for $(m-h)p>N>(m-h-1)p$ with $0<\l\le m-h-{N/p}$;
\item "(v)" $C^{h,\l}(Q)$ for $N=(m-h-1)p$ with 
$0<\l< 1$.
\endroster

Then, if $\T^{\a}_{m,k,p}({\Cal A})>0$, the following inequality is valid:

$$
||u||_T\le C_0\left(||\n^{k+1}u||_{L^{p_0}(Q)}
+C||\n^mu||_{L^p(Q)}\right).
\leqno(2.2)
$$

Here $C_0$ a (maybe large) constant, $C_0=C_0(N,m,p,p_0,\a)$, and
the constant $C$ can be chosen as

$$
c\cdot\T^{\a}_{m,k,p}({\Cal A})^{-{1\over p}}.
$$

Here $c$ is a constant $c=c(N,m,p,p_0,\a)$.
\par

On the other hand assume that the inequality {\rm (2.3)} is true for all
$u\in{\Cal A}$, with $T=L^q(Q)$ and $q$ as above. 
Furthermore $p_0$ is as $q$ above, but with $m$ replaced with $m-k-1$.
Then there are positive constants $c_1=c_1(N,m,p)$ and
$C_1=C_1(N,m,p,\a)$, such that

$$
||u||_{L^q(Q)}\le C_1\left(||\n^{k+1}u||_{L^{p_0}(Q)}
+||\n^mu||_{L^p(Q)}\right)
\leqno(2.3)
$$
or

$$
C\ge c_1\cdot\T^{\a}_{m,k,p}({\Cal A})^{-{1\over p}}.
$$

\endproclaim
\bs

First we give a lemma.
\bs

\proclaim{2.10. Lemma}
Let $Q$ be a fixed cube in \R. Let $p\ge 1$, let $q$ be as in Theorem
\text{\rm 2.9} and let $p_0>0$.
Let $u\in W^{m,p}(Q)$.
Then if 

$$
||u-\Pi_{m-1,m-1,p}u||_{L^p(Q)}
>{1\over 2}||\Pi_{m-1,m-1,p}u||_{L^p(Q)}
\leqno (2.4)
$$
or

$$
||\Pi_{m-1,k,p}u||_{L^p(Q)}
<\a||\Pi^{\text{co}}_{m-1,k,p}u||_{L^p(Q)},
\leqno (2.5)
$$
then

$$
||u||_{L^q(Q)}
\le A(||\n^{k+1}u||_{L^{p_0}(Q)}+||\n^mu||_{L^p(Q)}).
$$

In the first case $A=A(N,m,p)$ and in the second case $A=A(N,m,p,\a)$.
\endproclaim

\demo{Proof of Lemma}

First we assume (2.4).
We use the Sobolev inbedding theorem, 
the triangle inequality, (2.4) and a weak Poincar\'e inequality
to get

$$
\eq{
||u||_{L^q(Q)}
\le 
&A(||u||_{L^p(Q)}+||\n^mu||_{L^p(Q)})
\crd
\le
&A(||u-\Pi_{m-1,m-1,p}u||_{L^p(Q)}
+||\Pi_{m-1,m-1,p}u||_{L^p(Q)}
+||\n^mu||_{L^p(Q)})
\crd
<
&A(3||u-\Pi_{m-1,m-1,p}u||_{L^p(Q)}
+||\n^mu||_{L^p(Q)})
\crd
\le
&A||\n^mu||_{L^p(Q)}.
\cr}
$$

Next we assume (2.5).
We use equivalent norms in finite dimensional vector 
space, the triangle inequality and a weak Poincar\'e inequality
to get

$$
\eq{
||\Pi^{\text{co}}_{m-1,k,p}u||_{L^p(Q)}
\le
&A||\n^{k+1}\Pi^{\text{co}}_{m-1,k,p}u||_{L^p(Q)}
\crd
=&A||\n^{k+1}\Pi_{m-1,m-1,p}u||_{L^p(Q)}
\crd
=&A||\n^{k+1}(u-\Pi_{m-1,m-1,p}u)-\n^{k+1}u||_{L^p(Q)}
\crd
\le&
A(||\n^{k+1}(u-\Pi_{m-1,m-1,p}u)||_{L^p(Q)}
\crd
&+||\n^{k+1}u||_{L^p(Q)})
\crd
\le
&A(||\n^{k+1}u||_{L^p(Q)}+||\n^mu||_{L^p(Q)}).
\cr}\leqno(2.6)
$$

Now we have by the Sobolev imbedding theorem, 
the triangle inequality, a weak Poincar\'e inequality,
(2.5), (2.6) and the argument in Theorem 2.4 that

$$
\eq{
||u||_{L^q(Q)}
\le
&A(||u||_{L^p(Q)}
+||\n^mu||_{L^p(Q)})
\crd
\le
&A(||u-\Pi_{m-1,m-1,p}u||_{L^p(Q)}
+||\Pi^{\text{co}}_{m-1,k,p}u||_{L^p(Q)}
\crd
&+||\Pi_{m-1,k,p}u||_{L^p(Q)}
+||\n^mu||_{L^p(Q)})
\crd
\le
&A(||\Pi^{\text{co}}_{m-1,k,p}u||_{L^p(Q)}
+||\Pi_{m-1,k,p}u||_{L^p(Q)}
+||\n^mu||_{L^p(Q)})
\crd
\le
&(A+\a)||\Pi^{\text{co}}_{m-1,k,p}u||_{L^p(Q)}
+A||\n^mu||_{L^p(Q)}
\crd
\le
&A(1+\a)(||\n^{k+1}u||_{L^p(Q)}
+||\n^mu||_{L^p(Q)})
\crd
\le
&A(1+\a)(||\n^{k+1}u||_{L^{p_0}(Q)}
+||\n^mu||_{L^p(Q)}).
\cr}
$$
\vfill\eject

\demo{Proof of theorem}
We begin with the first part. 
By the same arguments used in Theorem 2.4 it is sufficient to prove
the inequality for $p_0=p$ and $T$ equal to $L^p(Q)$.
\par

According to Lemma 2.10 we can assume that

$$
||u-\Pi_{m-1,m-1,p}u||_{L^p(Q)}<{1\over 2}||\Pi_{m-1,m-1,p}u||_{L^p(Q)}
\leqno(2.7)
$$
and
$$
||\Pi_{m-1,k,p}u||_{L^p(Q)}\ge\a||\Pi^{\text{co}}_{m-1,k,p}u||_{L^p(Q)}.
\leqno(2.8)
$$

We have by the triangle inequality, (2.8), and (2.7) that

$$
\eq{
||u||_{L^p(Q)}
\le 
&||u-\Pi_{m-1,m-1,p}u||_{L^p(Q)}
+||\Pi^{\text{co}}_{m-1,k,p}u||_{L^p(Q)}
\crd
&+||\Pi_{m-1,k,p}u||_{L^p(Q)}
\crd
\le
&||u-\Pi_{m-1,m-1,p}u||_{L^p(Q)}
+\left(1+{1\over \a}\right)||\Pi_{m-1,k,p}u||_{L^p(Q)}
\crd
\le
&
\left({3\over 2}+{1\over \a}\right)||\Pi_{m-1,k,p}u||_{L^p(Q)}.
\cr}
$$
By (2.7) and (2.8) $u\in T^{\a}_{m,k,p,{\Cal A}}$.
Now by the definition of the polynomial capacity $\T^{\a}_{m,k,p}$
we have 

$$
||\Pi_{m-1,k,p}u||_{L^p(Q)}
\le
{||\n^mu||_{L^p(Q)}\over
\inf_{v\in T^{\a}_{m,k,p,{\Cal A}}}
{\dsize {||\n^mv||_{L^p(Q)}\over
||\Pi_{m-1,k,p}v||_{L^p(Q)}}}}
=
{||\n^mu||_{L^p(Q)}\over \T^{\a}_{m,k,p}({\Cal A})}.
$$

This ends the proof of the first part.
\bs

Next we prove the last part of the theorem. 
We can assume that $C$ is large. We assume that 

$$
||u||_{L^q(Q)}
\le
C_0(||\n^{k+1}u||_{L^{p_0}(Q)}
+C||\n^mu||_{L^p(Q)})
\leqno(2.9)
$$
for all $u\in{\Cal A}$. 
We have by equivalent norms (and quasinorms)
in finite dimensions, triangle inequalities, (2.9), the Sobolev imbedding
theorem and weak Poincar\'e inequalities that

$$
\eq{
\phantom{<}&||\Pi_{m-1,k,p}u||_{L^p(Q)}
\crd
\le
&A||\Pi_{m-1,k,p}u||_{L^q(Q)}
\crd
\le 
&A(||u-\Pi_{m-1,k,p}u||_{L^q(Q)}
+||u||_{L^q(Q)})
\crd
\le
&A(||u-\Pi_{m-1,k,p}u||_{L^q(Q)}
+||\n^{k+1}u||_{L^{p_0}(Q)}
+C||\n^mu||_{L^p(Q)})
\crd
\le
&A(||u-\Pi_{m-1,k,p}u||_{L^q(Q)}
+||\n^{k+1}u||_{L^p(Q)}
+C||\n^mu||_{L^p(Q)})
\crd
\le
&A(||u-\Pi_{m-1,k,p}u||_{L^p(Q)}
+||\n^{k+1}u||_{L^p(Q)}
+C||\n^mu||_{L^p(Q)})
\crd
\le
&A(||u-\Pi_{m-1,m-1,p}u||_{L^p(Q)}
+||\Pi^{\text{co}}_{m-1,k,p}u||_{L^p(Q)}
+||\n^{k+1}u||_{L^p(Q)}
\crd
\phantom {\le}
&+C||\n^mu||_{L^p(Q)})
\crd
\le
&A\left(
\Pi^{\text{co}}_{m-1,k,p}u||_{L^p(Q)}
+||\n^{k+1}u||_{L^p(Q)}
+C||\n^mu||_{L^p(Q)}
\right)
\crd
\le
&A(||\Pi^{\text{co}}_{m-1,k,p}u||_{L^p(Q)}
+||\n^{k+1}(u-\Pi_{m-1,m-1,p}u)||_{L^p(Q)}
\crd
\phantom{\le}
&+||\n^{k+1}\Pi^{\text{co}}_{m-1,k,p}u||_{L^p(Q)}
+C||\n^mu||_{L^p(Q)})
\crd
\le
&A(||\Pi^{\text{co}}_{m-1,k,p}u||_{L^p(Q)}
+C||\n^mu||_{L^p(Q)}).
\cr}
\leqno(2.10)
$$

It is enough to consider those $u$ satisfying (2.7) and (2.8)
by Lemma 2.10, i.e. we consider only the test functions for 
$\T^{\a}_{m,k,p}$.
Hence by (2.8) and (2.10) we have, taking a $u$ giving approximately
the value of $\T^{\a}_{m,k,p}({\Cal A})$, that

$$
\left(1-{A\over \a}\right)
{||\n^mu||_{L^p(Q)}\over \T^{\a}_{m,k,p}({\Cal A})^{1\over p}}
\le
A\cdot C||\n^mu||_{L^p(Q)}.
$$

End of proof.
\bs

Now we give the definition of the standard Sobolev condenser capacity
related to the cubes $Q$ and $2Q$ in \R\ and Sobolev space $W^{m,p}$.
\bs

\no{\smc 2.11. Definition.} Let $Q$ be a fixed cube in \R. Let $p\ge 1$,
Let $K$ be a compact set in the closure of $Q$ then the consenser
capacity

$$
C_{m,p}(K)=
\inf_{\{\varphi\in C^{\infty}_0(2Q):\varphi\ge 1_K\}}
\{||\n^m\varphi||^p_{L^p(2Q)}\}
$$
\bs

We will give several estimates of the polynomial capacities, i.e. of
the $\G$ and $\T$ types, together with ordinary capacities $C_{m,p}$.
\vfill\eject

We begin with a comparison between $\G_{m,k,p}$ and $\T^{\a}_{m,k,p}$.
\bs

\proclaim{2.12. Corollary of Theorems 2.4 and 2.9}

Let $p\ge 1$, let $Q$
be a fixed cube in \R\ and let ${\Cal A}\subset W^{m,p}(Q)$.
Then for $\G_{m,k,p}({\Cal A})$ and $\G_{m,m-1,p}({\Cal A})$ small
respectively

$$
\T^{\a}_{m,k,p}({\Cal A})\lesssim \G_{m,k,p}({\Cal A})
$$
and

$$
\T^{\a}_{m,m-1,p}({\Cal A})\sim \G_{m,m-1,p}({\Cal A}).
$$

The equivalences are independent of $\Cal A$.
\endproclaim

\demo{Proof} Compare Theorems 2.4 and 2.9.
\bs

We give the definition of another capacity similar to the condenser
capacity $C_{m,p}$.
\bs

\demo{2.13. Definition}
Let $p\ge 1$ and let $Q$ be an open cube in \R.
Let $K$ be a closed set in $\bar Q$.
Define
$$
C^{\sharp}_{m,p}(K)=\inf_{u\in S^{\sharp}_K}||\n^mu||_{L^p(Q)},
$$
where 
$$
S^{\sharp}_K=\{u\in C^{\infty}_0(2Q):u= 1\text{ in a neighbourhood of }K\}.
$$

The following result is due to Maz'ya, see [MAZ], and
D.R. Adams and J.C. Polking, see [A-P].

In fact it follows from Spectral Synthesis in Sobolev space, 
with $=$ instead of $\sim$ in the formulation. 

It is very much easier to make the direct proof.
\bs

\proclaim{2.14. Theorem} Let $p>1$, let $Q$ be an open cube in \R
and let $K$ be a closed set in
$\bar Q$ in \R. Then

$$
C^{\sharp}_{m,p}(K)\sim C_{m,p}(K).
$$

The equivalence is independent of $K$.
\endproclaim
\bs

Now we give some more background on the properties of Sobolev spaces.
\bs

The $L^p$ functions are defined only almost everywhere with respect to
Lebesgue measure. Lebesgue measure actually according to the definition
of capacities is equal to $C_{0,p}$ for any $p$. Hence it is quite
natural that higher order capacities takes the place of Lebesgue
measure as the order of Sobolev space gets higher. Zero order Sobolev
space is $L^p$ space. We have defined condenser capacities and rightly
there should be whole space definitions of the capacities to suit
the purpose. However we save the reader this definition and instead leave it
to he/she to use these condenser ones anyway, since it is possible.
\bs

There are some theorems that capitalizes on this connection.
\bs

The functions in $W^{m,p}$ can be defined quasi-everywhere (q.e.) i.e.
up to a set of $C_{m,p}$-capacity zero. How this works is made clear 
in the following definition and theorem, see [M-H].
\bs

\no{\smc 2.15. Definition.} A function $u\in W^{m,p}$ is called 
quasicontinuous if for each $\e>0$
there exists an open set $\omega$ such that 
$C_{m,p}(\omega)<\e$ and $u|_{\omega^c}$ is
continuous.
\bs

\proclaim
{2.16. Theorem} Let $p>1$. If $u\in W^{m,p}$ there exists a quasicontinuous
function $\bar u$ such that $u$ and $\bar u$ coincides almost everywhere.

If $\bar u$ and $\bar{\bar u}$ are quasicontinuous and coincides a.e.
then they coincides q.e.
\endproclaim

Now it is possible to define subspaces of $W^{m,p}$ with different 
trace properties properly, i.e. we assume $W^{m,p}$ to consist of
quasicontinuous functions only.
\bs

It is well-known that it does not change the infimum 
in Definition 2.11 and in Definition 2.13)
to take as test functions (or analogously)

$$
S^{\p}_K=\{u\in W^{m,p}_0(2Q):u|_K\ge 1\ \text{q.e.}\}
$$
instead of $S_K$. 
\bs

\no{\smc 2.17. Definition.}
Let $1<p$, let $0\le r\le m-1$
and let $\O$ be an open set in \R. Then define
$$
W^{m,p}_{0,r}(\O)=\{u\in W^{m,p}:D^{\a}u|_{\O^c}=0\ \text{q.e.}
\ |\a|\le r\}.
$$
(Here q.e. means: Except for a set $E_{\a}$, such that 
$C_{m-|\a|,p}(E_{\a})=0$.)
\bs

\proclaim{2.18. Proposition} Let $p\ge 1$, let $Q$ be a 
fixed closed cube in \R, let $K$ be a closed set in $Q$ and let 
$$
{\Cal A}=W^{m,p}_0(K^c)|_Q
$$ 
or 
$$
{\Cal A}=W^{m,p}_{0,s}(K^c)|_Q.
$$
Let $\g$ be multiindex, with $0<|\g|\le k$ {\rm(}and in the second 
case $|\g|\le s${\rm)}. Then

$$
\G_{m,k,p}({\Cal A})
\gtrsim 
\G_{m-|\g|,k-|\g|,p}(D^{\g}{\Cal A}).
$$

The equivalence is independent of $K$.

\endproclaim
\bs

We prove an analogous result for the capacity $\T^{\a}_{m,k,p}$
stated for general sets of Sobolev space functions.
\vfill\eject

\proclaim{2.19. Proposition} Let $Q$ be a fixed cube in \R.
Let $p\ge 1$, let $\Cal A\subset W^{m,p}(Q)$
let $\a$ be large, positive and independent of $\Cal A$.
Then, for $s$, $0<s\le k$, there is a multiindex $\g$, with $|\g|=s$,
there are $\a'$ and multiindex $\delta$, with $\delta\le \g$,
such that

$$
\aligned 
{}\\
\T^{\a}_{m,k,p}({\Cal A})\gtrsim
\\
{}
\endaligned
\aligned
\ &\T^{\a'}_{m-s,k-s,p}(D^{\g}{\Cal A})
\\
\ &\text{or}
\\
\ &\T^{\a'}_{m-|\delta|,0,p}(D^{\delta}{\Cal A}).
\endaligned
$$
Here the equivalence depends on $N$, $m$, $p$ and $\a$.

\endproclaim

\demo{Proof} 
Either $\T^{\a}_{m,k,p}({\Cal A})=1$, and then 
we are done, 
or $\T^{\a}_{m,k,p}({\Cal A})<1$.
In this case
it is sufficient to prove the estimate below for a multiindex $\g$,
with $|\g|=1$, since the general case follows by induction,

$$
\aligned
{}
\\
A\cdot\T^{\a}_{m,k,p}(\Cal A)\ge 
\\
{}
\endaligned
\aligned
\ &\T^{\a''}_{m-|\g|,k-|\g|,p}(D^{\g}{\Cal A})
\\
\ &\text{or}
\\
\ &\T^{\a''}_{m,0,p}(\Cal A).
\endaligned
\leqno(2.11)
$$
Here $\a''$ is some positive constant independent of $\Cal A$.
\par

Take a test function $u$ for 
$\T^{\a}_{m,k,p}({\Cal A})$, which gives approximately the right
value of this capacity.
We want to estimate the quantity

$$
{||\n^mu||_{L^p(Q)}\over ||\Pi_{m-1,k,p}u||_{L^p(Q)}}.
\leqno(2.12)
$$

Assume that 

$$
||\Pi_{m-1,0,p}u||_{L^p(Q)}> ||\Pi^{\p}_{m-1,k,p}u||_{L^p(Q)}.
\leqno(2.13)
$$

Then by equivalent norms in finite dimensional space and
(2.13) we have that

$$
{||\n^mu||_{L^p(Q)}\over ||\Pi_{m-1,k,p}u||_{L^p(Q)}}
\sim
{||\n^mu||_{L^p(Q)}\over
||\Pi'_{m-1,k,p}u||+||\Pi_{m-1,0,p}u||_{L^p(Q)}}
\sim
{||\n^mu||_{L^p(Q)}\over ||\Pi_{m-1,0,p}u||_{L^p(Q)}}.
$$
\bs

It is clear by (2.13) that $u\in T^{\a''}_{m,0,p,{\Cal A}}$
for some constant $\a''$ independent of $\Cal A$.
Hence it follows that in the case considered
\bs

$$
\T^{\a}_{m,k,p}({\Cal A})
\sim   
\T^{\a''}_{m,0,p}({\Cal A}).
$$
\vfill\eject

On the other hand assume that 

$$
||\Pi_{m-1,0,p}u||_{L^p(Q)}\le ||\Pi^{\p}_{m-1,k,p}u||_{L^p(Q)}.
\leqno(2.14)
$$

We make the following calculation to be used later.
We have by the triangle inequality and weak 
Poincar\'e inequalities that
$$
\eq{
&|\ ||\Pi_{m-2,m-2,p}D^{\g}u||_{L^p(Q)}
-||D^{\g}\Pi_{m-1,m-1,p}u||_{L^p(Q)}\ |
\crd
\le 
&||\Pi_{m-2,m-2,p}D^{\g}u-D^{\g}\Pi_{m-1,m-1,p}u||_{L^p(Q)}
\crd
=
&||(D^{\g}u-\Pi_{m-2,m-2,p}D^{\g}u)-(D^{\g}u-D^{\g}\Pi_{m-1,m-1,p}u)||_{L^p(Q)}
\crd
\le&
||D^{\g}u-\Pi_{m-2,m-2,p}D^{\g}u||_{L^p(Q)}+
||D^{\g}u-D^{\g}\Pi_{m-1,m-1,p}u||_{L^p(Q)}
\crd
\le
&A(||\n^{m-1}D^{\g}u||_{L^p(Q)}
+||\n^mu||_{L^p(Q)})
\crd
\le&
A_0||\n^mu||_{L^p(Q)}.
\cr}
\leqno(2.15)
$$

The expression (2.12) can be rewritten as

$$
{||\n^mu||_{L^p(Q)}\over ||\Pi_{m-2,k-1,p}D^{\g}u||_{L^p(Q)}}
\cdot
{||\Pi_{m-2,k-1,p}D^{\g}u||_{L^p(Q)}
\over
||\Pi_{m-1,k,p}u||_{L^p(Q)}}.
\leqno(2.16)
$$

The procedure is first to estimate the first factor from below, with 
$\T^{\a''}_{m,k,p}({\Cal A})$ and then the second with $1$.
The first step in this is to show that, with the right choice
of $\g$ and $\a''$,

$$
\eq{
&D^{\g}u\in T^{\a''}_{m-1,k-1,p,D^{\g}{\Cal A}}
\crd
&\text{or}
\crd
&\T^{\a}_{m,k,p}({\Cal A})\sim 1.
\cr}
\leqno(2.17)
$$

In the second case there is nothing to prove.
\par

First we check that

$$
||\Pi_{m-2,k-1,p}D^{\g}u||_{L^p(Q)}
\ge 
\a''||\Pi^{\text{co}}_{m-2,k-1,p}D^{\g}u||_{L^p(Q)}.
\leqno(2.18)
$$

We know that 

$$
||\Pi_{m-1,k,p}u||_{L^p(Q)}
\ge 
\a||\Pi^{\text{co}}_{m-1,k,p}u||_{L^p(Q)}.
\leqno(2.19)
$$

The idea is to compare the LHS of (2.18) and (2.19) as well as
the RHS of (2.18) and (2.19)
respectively. To begin with, by (2.14) and equivalent norms in finite
dimensional space and with the right choice of $\g$ we have that

$$
||\Pi_{m-1,k,p}u||_{L^p(Q)}\sim
||\Pi'_{m-1,k,p}u||_{L^p(Q)}
\sim
||D^{\g}\Pi_{m-1,k,p}u||_{L^p(Q)}.
\leqno(2.20)
$$

Now assume 

$$
||\Pi_{m-2,k-1,p}D^{\g}u||_{L^p(Q)}
\le
2A_0||\n^mu||_{L^p(Q)}.
\leqno(2.21)
$$

Then by (2.15), (2.21) and the triangle inequality we have that

$$
||D^{\g}\Pi_{m-1,k,p}u||_{L^p(Q)}
\le
3A_0 ||\n^mu||_{L^p(Q)}.\leqno(2.22)
$$
But by (2.22) and (2.20) we have that

$$
||\Pi_{m-1,k,p}u||_{L^p(Q)}
\lesssim
||\n^mu||_{L^p(Q)}
$$
and this implies that $\T^{\a}_{m,k,p}({\Cal A})\sim 1$,
i.e. this case is done.

If the inequality opposite to (2.21) is true then by (2.15)

$$
||\Pi_{m-2,k-1,p}D^{\g}u||_{L^p(Q)}
\sim
||D^{\g}\Pi_{m-1,k,p}u||_{L^p(Q)},\leqno(2.23)
$$
and by (2.20) and (2.23) we get the desired estimate.
\bs

We do a similar
estimate of the RHS of (2.18).
We begin a chain of estimates by the observation that, by 
equivalent norms in finite dimensional space,

$$
||\Pi^{\text{co}}_{m-1,k,p}u||_{L^p(Q)}
\gtrsim
||D^{\g}\Pi^{\text{co}}_{m-1,k,p}u||_{L^p(Q)}.
\leqno(2.24)
$$

Assume that 

$$
||\Pi^{\text{co}}_{m-2,k-1,p}D^{\g}u||_{L^p(Q)}
\ge
2A_0||\n^mu||_{L^p(Q)}.
\leqno(2.25)
$$

Then by (2.15), and (2.25)

$$
||\Pi^{\text{co}}_{m-2,k-1,p}D^{\g}u||_{L^p(Q)}
\sim
||D^{\g}\Pi^{\text{co}}_{m-1,k,p}u||_{L^p(Q)}.
\leqno(2.26)
$$

The desired estimate follows from (2.24) and (2.26).
\par
On the other hand assume the inequality opposite to (2.25);

$$
2A_0 ||\n^mu||_{L^p(Q)}
>
||\Pi^{\text{co}}_{m-2,k-1,p}D^{\g}u||_{L^p(Q)}.
\leqno(2.27)
$$

We assume that we have the opposite to what we want, i.e.

$$
||\Pi_{m-2,k-1,p}D^{\g}u||_{L^p(Q)}
<
\a''||\Pi^{\text{co}}_{m-2,k-1,p}D^{\g}u||_{L^p(Q)}.
\leqno(2.28)
$$

By (2.27), (2.28) and (2.15) we have that

$$
\eq{
3A_0||\n^mu||_{L^p(Q)}
\ge
&{1\over\a''}||\Pi_{m-2,k-1,p}D^{\g}u||_{L^p(Q)}
\crd
\ge
&{1\over\a''}\left(||D^{\g}\Pi_{m-1,k,p}u||_{L^p(Q)}
-A_0||\n^mu||_{L^p(Q)}\right).
\cr}
\leqno(2.29)
$$

By (2.29) and (2.20) we have, if $\a''$ is large enough, that

$$
||\n^mu||_{L^p(Q)}
\gtrsim
||\Pi_{m-1,k,p}u||_{L^p(Q)}
$$

and this proves that $\T^{\a}_{m,k,p}({\Cal A})\sim 1$ in this case.
\bs
In order to prove (2.17)  we also
have to study

$$
||D^{\g}u-\Pi_{m-2,m-2,p}D^{\g}u||_{L^p(Q)}
\le
{1\over 2}||\Pi_{m-2,m-2,p}D^{\g}u||_{L^p(Q)}.
\leqno(2.30)
$$

Assume the opposite of (2.30). 
By Lemma 2.10. and Theorem 2.9 we have that 
$\T^{\a}_{m,k,p}({\Cal A})\sim 1$.
All this proves (2.17).
Thus we have proved that either the first factor of (2.16) is less
than or equal to $\T^{\a''}_{m-1,k-1,p}(D^{\g}{\Cal A})$ or that
$\T^{\a}_{m,k,p}({\Cal A})\sim 1$.

It remains to estimate the second factor in the RHS of (2.16):

$$
{||\Pi_{m-2,k-1,p}D^{\g}u||_{L^p(Q)}
\over
||\Pi_{m-1,k,p}u||_{L^p(Q)}}.
\leqno(2.31)
$$

First we assume (2.21). In this case we already are done since it
implies that $\T^{\a}_{m,k,p}({\Cal A})\sim 1$.
Assume the opposite to (2.21) then by 
(2.23) and (2.20) we get the desired result.
\par

End of proof.
\bs

\no{\smc 2.20. Remark.} We see here that there often appears
conditions on $\a$ being large enough. It is seen from Theorem
2.9 that $\T^{\a}_{m,k,p}({\Cal A})\sim \T^{\b}_{m,k,p}({\Cal A})$,
when $\a,\b\ge a$ for some $a=a(N,m,p)$. The equivalence is independent
of $\Cal A$.
\vfill\eject

The $\G$-capacities can be compared with the capacity $C_{m,p}$.

\proclaim{2.21. Proposition}
Let $p> 1$ and let $K$ be a closed set in the fixed 
closed cube $Q$ in \R. Then

$$
\G_{m,0,p}(W^{m,p}_{0,0}(K^c)|_Q)
\sim \G_{m,0,p}(W^{m,p}_0(K^c)|_Q)
\sim C_{m,p}(K).
$$
\endproclaim
\bs

\demo{Proof} Follows from the definition of the $\G$-capacity,
Theorem 2.14 and [MEY2; Th. 2.5].
\bs

We will now show a similar estimate for $\T^{\a}_{m,k,p}$.
\bs

\proclaim{2.22. Proposition} Let $p\ge 1$, let $Q$ be a fixed open cube in
\R\ and let $K$ be a closed set
in $\bar Q$. If $\a>0$ is large, positive and 

$$
{\Cal A}=W^{m,p}_0(K^c)|_Q
$$
or 
$$
{\Cal A}=W^{m,p}_{0,0}(K^c)|_Q,
$$

then

$$
A\cdot\T^{\a}_{m,0,p}({\Cal A})
\ge 
C_{m,p}(K).
$$

Here $A$ depends only on $N$, $m$, $p$ and  $\a$.
\endproclaim
\bs

\demo{Proof} We can assume $\T^{\a}_{m,0,p}({\Cal A})<1$, since
otherwise the assertion is trivial because the capacities $C_{m,p}(K,2Q)$
have a fixed upper bound. Hence $\T^{\a}_{m,0,p}({\Cal A})$
is determined by test functions belonging to
$T^{\a}_{m,0,p,{\Cal A}}$. Let $u\in T^{\a}_{m,0,p,{\Cal A}}$.
By a trivial estimate, the Hestenes construction and a weak Poincar\'e
inequality we get

$$
\eq{
||\n^mu||_{L^p(Q)}
&\le
||\n^m(u-\Pi_{m-1,m-1,p}u)^*||_{L^p}
\crd
&\le
A\sum_{i=0}^{m-1}||\n^i(u-\Pi_{m-1,m-1,p}u)||_{L^p(Q)}
\crd
&\le
A||\n^mu||_{L^p(Q)}.
\cr}
$$

Hence 

$$
{||\n^mu||_{L^p(Q)}
\over
 ||\Pi_{m-1,0,p}u||_{L^p(Q)}}
\sim
{||\n^m(u-\Pi_{m-1,m-1,p}u)^*||_{L^p}
\over
||\Pi_{m-1,0,p}u||_{L^p(Q)}}.
$$

Now by equivalent norms in finite dimensional space we have that

$$
||\Pi^{\text{co}}_{m-1,0,p}u||_{L^p(Q)}
\sim
||\Pi^{\text{co}}_{m-1,0,p}u||_{L^{\infty}(Q)}.
$$
\vfill\eject

Trivially 

$$
||\Pi_{m-1,0,p}u||_{L^p(Q)}
=\text{const}\cdot||\Pi_{m-1,0,p}u||_{L^{\infty}}.
$$

It follows that if $\a$ is large we have that

$$
v={-(u-\Pi_{m-1,m-1,p}u)^*\over \Pi_{m-1,0,p}u}
\biggl|_K\ge {1\over 2}\text{\ q.e.}
$$

\no
say. This implies that $2v$ can be used as test function for 
$C_{m,p}(K,2Q)$. It follows that if $\a$ is large  enough,
then

$$
C_{m,p}(K,2Q)\le A\cdot\T^{\a}_{m,0,p}({\Cal A}).
$$

End of proof.
\bs

Summing up we have a consequence of Theorem 2.9 and 
of the properties of the capacity $\T^{\a}_{m,k,p}$.

\proclaim
{2.23. Corollary} Let $p\ge 1$. Let $Q$ be an open cube in \R and
let $K$ be a closed subset of $\bar Q$. Let $u\in W^{m,p}_{0,k}(K^c)|_Q$.
Then there is a constant $C=C(N,m,p)$ such that 

$$
||u||_{L^p(Q)}
\le 
C||\n^{k+1}u||_{L^p(Q)}
+{C\over C_{m-k,p}(K,2Q)} ||\n^mu||_{L^p(Q)}.
$$
\endproclaim

\demo{Proof} Follows from Theorem 2.22, Proposition 2.18 and 
Proposition 2.22.
\bs

\no{\smc 2.24. Remark.} Corollary 2.23 was proved by Hedberg for
$p>1$, see [HED]. That result was the inspiration to treat this
Poincar\'e inequality in a complete kind of way.
\bs

The polynomial capacity, $\G_{m,k,p}(W^{m,p}_0(K^c))|_Q$,
can sometimes get a better estimate if the boundary is sufficiently
irregular i.e. deviates from the algebraic surfaces of a certain
degree. Of the same ought to hold for $\T^{\a}_{m,k,p}$ as well but 
the proof has not been carried out.
\bs

\proclaim{2.25. Proposition}
Let $p\ge 1$, let $0\le k,s\le m-1$ and 
let $Q$ be an open cube in \R. 
Let $\{x_i\}$, a finite set of points in $\bar Q$, satisfy that
for any polynomial $P\in {\Cal P}_{k-l}$, where $0\le l\le s,k$,
at least one of the $x_i$ is at least at distance $d$ from
the zero set of $P$.
Let $K$ be a closed set in $\bar Q$. 
Then there is an $r=r(N,m,p,d)$ such that 

$$
C\cdot\G_{m,k,p}(W^{m,p}_0(K^c)|_Q)\ge \min_iC_{m-l,p}(K\cap B_r(x_i));
$$
$$
C\cdot\G_{m,k,p}(W^{m,p}_{0,s}(K^c)|_Q)\ge \min_iC_{m-l,p}(K\cap B_i).
$$
Here $C$ is independent of $K$.
\endproclaim
\bs

\demo{Proof}
We prove only the second estimate. The first one can be seen as a
consequence, but also can be proved in the same way.
First we use Proposition 2.17 to get
$$
A\cdot\G_{m,k,p}(W^{m,p}_{0,s}(K^c)|_Q)
\ge
\G_{m-l,k-l,p}(W^{m-l,p}_{0,0}(K^c)|_Q).
$$
\par

We can assume that every $P$ has $||P||_{L^p(Q)}=1$.
Now $\max_i|P(x_i)|$ constitutes a norm on ${\Cal P}_{k-l}$.
Hence by equivalent norms in finite dimensional space
${1\over a}\ge \max_i|P(x_i)|\ge a>0$ for a constant 
$a=a(N,m,p,d)$.
Furthermore $0\le ||\n P||_{L^{\infty}}\le b$ for a constant 
$b=b(N,m,p,d)$ also by equivalent norms in finite dimensions.
(Adjusted to the seminorm case of course.)
\par

Say $i_0$ gives the maximimum above.
Put $r={a\over 2b}$. Then 

$$
r||\n P||_{L^{\infty}(Q)}\le {a\over 2}\le {1\over 2}|P(x_{i_0})|
$$
and it then follows that $P(x)\ge {a\over 2}$ in $B_r(x_{i_0})$.
\par

Finally we note that if $P$ and $u$ determine approximately 
\pn
$\G_{m-l,k-l,p}(W^{m-l,p}_{0,0}(K^c)|_Q)$,
then we are done since $2u/a$ is a test function for the capacity
$C_{m-l, p}(K\cap B_{i_0})$, 
(which is a Sobolev space function instead of a $C^{\infty}$
function).
\par

End of proof.
\bs

The Proposition 2.25. should be compared with a result on sets with
the Markov property by Nystr\"om which has been proved later, see [NYS].
\bs

Here we more modestly apply it to a self similar set, which is special
case.
\bs

We show that a self simlar set, which is not subset of a hyperplane 
deviates uniformly with respect to dilations 
from algebraic surfaces of a fixed degree.
\bs

\proclaim {2.26. Proposition} Let $S$ be closed set in \R, which is
self similar and not a subset of a hyperplane.
Let $\Cal F$ be a Whitney decomposition of $S^c$.
Then 

$$
\inf_{Q\in {\Cal F}}
\G_{m-1,m-1,p}(W^{m,p}_{0,0}(\tilde R_Q\cap \tilde S))>0.
$$
\endproclaim

\demo{Proof} Let $P\in {\Cal P}_{m-1}$ and denote with 
$V_P$ the algebraic surface defined as the zero set of $P$.
Define 

$$
F(P)=
\inf_{Q\in {\Cal F}}\ \sup_{x\in \tilde S\cap \tilde R_Q}
\ \inf_{y\in V_P}\ d(x,y).
$$
According to Proposition 4.22 and the definition of self similarity
it is enough to prove that $F(P)$
is bounded from below by a positive constant.
The polynomials can be assumed to have norm 1. They form a finite
dimensional unit sphere in ${\Cal P}_{m-1}$ and
is compact. $F(P)$
varies continuously with the coefficients of the polynomial.
Hence by compactness the minimum is attained by some 
$P_0\in{\Cal P}_{m-1}$ and it is enough to show that 
$\tilde S\cap \tilde R_Q\not\in V_{P_0}$.
Since $\tilde S\cap \tilde R_Q$ 
is not a subset of a hyperplane it contains points forming 
the extremal points of a $N$-dimensional simplex
within the ball defining the self similarity. This simplex
has a scaled image close to every point of $\tilde S$ and if 
$\tilde S\cap \tilde R_Q\in V_{P_0}$
then $V_{P_0}$ has singular points in all $\tilde S\cap \tilde R_Q$.
The set of singular points forms an algebraic variety 
(or a union of finitely many), $V^(1)_{P_0}$
of lower dimension. But now $\tilde S\cap \tilde R_Q\in V^(1)_{P_0}$.
Proceeding inductively we get that $\tilde S\cap \tilde R_Q\in \emp$
A contradiction, hence the Hausdorff distance is greater than zero.
\par 

End of proof.
\bs
 
The following theorem shows that if we know that the set
${\Cal A}\subset W^{m,p}(Q)$ consists of nonnegative functions, 
then the $\T$-capacity sometimes can be estimated in better way than 
for similar subsets that do not have the nonnegative property.
\bs

\proclaim{2.27. Theorem}
Let $p\ge 1$, let $Q$ be  a fixed cube in \R\ and let 
${\Cal A}\subset W^{m,p}(Q)_+$.
Then for $\a$ large and some positive $\a'$

$$
\T^{\a}_{m,1,p}({\Cal A})\gtrsim\T^{\a'}_{m,0,p}({\Cal A}).
$$
The equivalence depends on $N$, $m$, $p$ and $\a$.
\bs

\endproclaim

\demo{Proof}
Suppose $u$ is a test function for the capacity,
$\T^{\a}_{m,1,p}({\Cal A})$, giving approximately the infimum.
First we show that if $\a$ is 
large enough, then

$$
||\Pi'_{m-1,1,p}u||_{L^p(Q)}
\le
A\left(||\Pi_{m-1,0,p}u||_{L^p(Q)}
+||\n^mu||_{L^p(Q)}\right).
\leqno(2.32)
$$

This follows from the following chain of inequalities.
Let $\text{neg}Q$ denote the part of $Q$, where $\Pi'_{m-1,1,p}u$
is negative and let $\text{pos}Q$ denote the part of $Q$, where $\Pi'_{m-1,1,p}u$
is positive. We have by the triangle inequality,
symmetry, nonnegativity of $u$, a weak Poincar\'e inequality and the
definition of the test functions, that
\vfill\eject

$$
\eq{
||\Pi'_{m-1,1,p}u||_{L^p(Q)}
\le
&||\Pi'_{m-1,1,p}u||_{L^p(\text{neg}Q)}
+||\Pi'_{m-1,1,p}u||_{L^p(\text{pos}Q)}
\crd
=
&2||-\Pi'_{m-1,1,p}u||_{L^p(\text{neg}Q)}
\crd
\le
&2||u-\Pi'_{m-1,1,p}u||_{L^p(\text{neg}Q)}
\crd
\le
&2||u-\Pi'_{m-1,1,p}u||_{L^p(Q)}
\crd
\le
&2||u-\Pi_{m-1,m-1,p}u||_{L^p(Q)}
+2||\Pi_{m-1,0,p}u||_{L^p(Q)}
\crd
&+2||\Pi^{\text{co}}_{m-1,1,p}u||_{L^p(Q)}
\crd
\le
&A||\n^mu||_{L^p(Q)}
+2||\Pi_{m-1,0,p}u||_{L^p(Q)}
+2||\Pi^{\text{co}}_{m-1,1,p}u||_{L^p(Q)}
\crd
\le
&A||\n^mu||_{L^p(Q)}
+2||\Pi_{m-1,0,p}u||_{L^p(Q)}
\crd
&+{2\over\a}||\Pi_{m-1,1,p}u||_{L^p(Q)}
\crd
\le
&A||\n^mu||_{L^p(Q)}
+\left(2+{2\over\a}\right)||\Pi_{m-1,0,p}u||_{L^p(Q)}
\crd
&+{2\over\a}||\Pi'_{m-1,1,p}u||_{L^p(Q)}.
\cr}
$$

Then choose $\a>2$ and (2.32) follows after rearrangement. 
\par
Now assume that

$$
||\n^mu||_{L^p(Q)}
\ge
||\Pi_{m-1,0,p}u||_{L^p(Q)}.
\leqno(2.33)
$$

Then by (2.32) we have that

$$
A||\n^mu||_{L^p(Q)}
\ge
||\Pi'_{m-1,1,p}u||_{L^p(Q)}
+||\Pi_{m-1,0,p}u||_{L^p(Q)}
\ge
||\Pi_{m-1,1,p}u||_{L^p(Q)},
$$
i.e. $\T^{\a}_{m,1,p}({\Cal A})\ge {1/A}$ and this case is done.
\par
Next assume

$$
||\n^mu||_{L^p(Q)}
<
||\Pi_{m-1,0,p}u||_{L^p(Q)}.
\leqno(2.34)
$$

Now from (2.32) it follows that

$$
A||\Pi_{m-1,0,p}u||_{L^p(Q)}
\ge
||\Pi_{m-1,1,p}u||_{L^p(Q)}
\leqno(2.35)
$$
and it follows that we can take $\a'=A(1+\a)$. Hence we have that 
$u\in T^{\a'}_{m,0,p,{\Cal A}}$.
Hence we get the theorem immediately from (2.35) and the definition of
$\T^{\a'}_{m,0,p}({\Cal A})$.
\par
End of proof.
\bs

We give a Poincar\'e type inequality for certain nonnegative functions.
It is a corollary of the earlier results.

{\bf 2.28. Theorem} Let $K$ be a closed subset of $\bar Q$ a
fixed cube in \R. Let $u\in W^{2,p}_{0,0}(K^c)$ with $u\ge 0$.
Then

$$
||u||_{L^p(Q)}
\le 
{A\over C_{2,p}(K)^{1\over p}}||\n^2u||_{L^p(Q)}.
\leqno(2.36)
$$
\bs

{\bf Proof:} Follows from Theorems 2.9, 2.27 and Proposition 2.22. 
\bs

The in a sense best constant.
\bs\bs

{\smc Section 3. On Spectral Synthesis in Sobolev Space}
\bs

To begin with we give one definition of Spectral synthesis in Sobolev
space. 
There is a dual formulation as well that makes the name adequate.
\bs

\no{\smc 3.1. Definition.} Spectral synthesis in Sobolev space is the fact:
Let $K$ be a closed set in \R and $m$, $p$ general, then 

$$
W^{m,p}_0(K^c)=W^{m,p}_{0,m-1}(K^c).
$$
\bs

Hedberg proved, see [HED] or better [A-H], that Spectral synthesis holds for every
closed set $K$ for $m\ge 1$, $p>1$. A result by Wolff was essential
in order to bring the possible values for $p$ all the way down.
\bs

The case $m=1$, $p\ge 1$ is much simpler and was proved earlier by Hedberg.
\bs

Hedberg actually proved a somewhat stronger formulation that the
result (from RHS to LHS) can be achieved by the use of a sequence a 
multipliers to generate a Cauchy sequence.

Netrusov, see [A-H], has proved this theorem in such a way that it
generalizes to (many) more function spaces as well.

Both these proofs are very long. 
\bs

The main interest lies in the original formulation however, but the 
multiplier formulation has had applications.
\bs

The Sobolev spectral synthesis question has a long history and it 
originates from Sobolev. He proved that if $K$ above is a Lipschitz
manifold then the conclusion follows. This may serve as an indication
of the progress made by Hedberg showing this for any closed $K$.
\bs

Proving Spectral synthesis is mostly thought of as proving the two 
inclusions that make up to the identity of the subspaces.

As it happens it is simple to prove 

$$
W^{m,p}_0(K^c)\subset W^{m,p}_{0,m-1}(K^c).
$$

We refer to Hedberg, see [HED2].
\bs

Now the remaining objective is to prove

$$
W^{m,p}_{0,m-1}(K^c)\subset W^{m,p}_0(K^c).
\leqno(3.1)
$$
\bs

Our idea here is to introduce a new technique to show the inclusion
(3.1) and to show it without multipliers.

However we fall short of this goal, but show a very interesting non-trivial
equivalent formulation. 

This formulation actually is powerful enough to give a corresponding 
theorem for nonnegative functions with $m=2$. This result is new.

This proof is the only proof with a $m=2$ situation as above with a short
proof. Otherwise only very long and complicated proofs are at hand
for $m>1$ as said before.
In this proof technique the polynomial capacities play an essential
role. 
\bs

Next we turn to a description of the technique itself.
\bs

It is enough to show that for any 
$u\in W^{m,p}_{0,m-1}(K^c)$ 
it is 
possible to approximate $u$ arbitrarily well by functions in $W^{m,p}_0(K^c)$,
since these subspaces are closed and that implies that we are done.
\bs

To simplify the proceedure decompose \R\ into equal sized 
disjoint cubes $Q_n$ and make a smooth partition of unity $\{\Phi_n\}$.
Hence it is obviously enough to consider $K$ compact.
\bs

One underlying important fact will be made used, namely the
fact that the two subspaces considered are equal when restriction is
taken to the closure of an open set. The closure not intersecting $K$. 

Now the idea is to piecewise interchange $u$ into $u'$
in a neighbourhood of $K$ only, in such a way that $u-u'$ has the
Sobolev norm is arbitrarily small, $u'\in W^{m,p}_0(K^c)$ and that way
get a proof.
\bs

Since $K$ can be taken compact then the function $u$ can be 
assumed to have bounded support. By 
a standard Poincar\'e inequality it follows that the Sobolev norm is
equivalent to the norm $||\n^mu||_{L^p}$ \ -- the norm used below.
\bs

Let $Q$ be an arbitrary cube with side $\delta$ that intersects $K$.

Suppose we are able to redefine

$$
u\in W^{m,p}_{0,m-1}(K^c)
$$ 

on $Q$ and get $u_1$ such that 

$$
u_1\in W^{m,p}_{0,m-1}(K^c),
$$

$$
u_1\in W^{m,p}_0(K^c)|_{{1\over 2}Q}
$$

and that 

$$
||\n^mu_1||_{L^p(Q)}\le A_Q||\n^mu||_{L^p(Q)}
\leqno(3.2)
$$

Next suppose that it is possible to iterate this for more cubes
in such way that if $u_j$ say has a certain $O_j$ set defined by this
proceedure such that

$$
u_j|_{O_j}\in W^{m,p}_0(K^c)|_{O_j},
$$
\bs

(*) \ {\bf Assume} that redefining iteratively gives $O_j\subset O_{j+1}$.
\bs

This means that the redefining process never destroys what
accomplished at an earlier stage.
\bs

After say $r$ iterations of this procedure $K$ is covered by $O_r$ and
$u_r$ is the wanted $u'$ approximant. 
\bs

(*) \ {\bf Assume} that the cubes covers $K$ at most $G$ times.
\bs

(*) \ {\bf Assume} that $A_Q=A_0$ for all cubes $Q$. 
\bs

Then we have the following estimate

$$
\eq{||\n^m(u_r-u)||_{L^p(\text{\bf R}^N)}^p
\le&
\sum_{j=1}^r
||\n^m(u_r-u)||_{L^p(Q_j)}^p
\crd
\le&
\sum_{j=1}^r
A(||\n^mu_r||_{L^p(Q_j)}^p+||\n^mu||_{L^p(Q_j)}^p),
\cr}
$$

but by the assumption on $A_Q$ and the intersection number $G$ it
follows that

$$
\eq{
\sum_{j=1}^r
A(||\n^mu_r||_{L^p(Q_j)}^p+||\n^mu||_{L^p(Q_j)}^p)
\le&
A\cdot (1+A_0)
\sum_{j=1}^r
||\n^mu||_{L^p(Q_j)}^p
\crd
\le&
A\cdot G\cdot (1+A_0)
||\n^mu||_{L^p(O_r)}^p.
\cr}
$$

Now $O_r$ shrinks with $\delta$ tending to zero so that

$$
||\n^mu||_{L^p(O_r)}^p
$$

tends to

$$
||\n^mu||_{L^p(K)}^p.
$$

However it can be shown (not difficult) that the truncation property
of Sobolev space of order one has the consequence that

$$
W^{m,p}_{0,m-1}(K^c)
=
W^{m,p}_{0,m}(K^c).
$$

Hence 

$$
||\n^mu||_{L^p(K)}^p=0,
$$

\no
i.e the approximation can be made arbitrarily good. 
\bs

This way we see the structure of the proof, also some parts of the
proof are done.
\bs

\proclaim{3.2 Fact} Let $K$ be a closed  set in \R. The
subspaces $W^{m,p}_0(K^c)$ and $W^{m,p}_{0,m-1}(K^c)$ are equal on
closure of open sets such that the closure does not intersect $K$.
\endproclaim

This follows by aproper use of a cut-off function.
\bs

Now we will concentrate on the task of redefining $u$ inside a given cube
$Q$ as described in the part of the proof above.
\bs

Let \R\ be tesselated by a lattice of cubes $\{Q'\}$ each of unit
size. In the standard way we construct functions 
$\phi_{Q'}\in C^{\infty}_0(Q')$ with $\phi|_{{1\over 2}Q'}=1$,
such that the $\phi_{Q'}$ are all equal except for a translation
allowed by the lattice.
Furthermore the set $\{\phi_{Q'}\}$ can be constructed to be a
partition of unity.

Next this lattice and functions are dilated so that the cube
sidelength is $\delta$. We denote the cubes $\{Q\}$ etc.
\bs

To begin with we only discuss the cubes $Q$ that intersect $K$.
\bs

Let $Q$ be such a cube. We want to redefine $u$ on $Q$ as indicated
previously. To begin with we do this first iteration discussed before.

Now some notation needeed.
\bs

For $f\in W^{m,p}(Q)$ define $\Pi f$ as the $P\in {\Cal P}_{m-1}$
that gives a minimum in

$$
\inf\ \{\sum^{m-1}_{i=0}l(Q)^{-(m-i)}||\n^i(f-P)||_{L^p(Q)}:
       P\in {\Cal P}_{m-1}\},
$$
where $l(Q)$ denotes the length of the side of the cube $Q$.
\bs

Say $v$ gives minimum in 

$$
\inf\ \{||\n^mf||_{L^p(Q)}:f\in W^{m,p}_{0,m-1}(K^c)|_Q,\ \Pi f=\Pi u\}
$$
and $w$ gives minimum in 

$$
\inf\ \{||\n^mf||_{L^p(Q)}:f\in W^{m,p}_0(K^c)|_Q,\ \Pi f=\Pi u\}.
$$

Then 

$$
||\n^mw||_{L^p(Q)}\le A(m,p,N)||\n^mv||_{L^p(Q)}\leqno(3.2)
$$

implies, first of all that $(1-\phi_Q)u+\phi_Qw$ has the right boundary
properties for the wanted redefinition. In order to show that the norm
properties also are right there is some calculation to do as follows

$$
\eq{
&||\n^m((1-\phi_Q)u+\phi_Qw)||_{L^p(Q)}
=||\n^m((1-\phi_Q)(u-\Pi u)
+\phi_Q(w-\Pi u))||_{L^p(Q)}
\crd
&\le
||\n^m((1-\phi_Q)(u-\Pi u))||_{L^p(Q)}
+||\n^m(\phi_Q(w-\Pi u))||_{L^p(Q)}
\crd
&\le
A\sum^m_{i=0}l(Q)^{-(m-i)}(||\n^i(u-\Pi u)||_{L^p(Q)}
+||\n^i(w-\Pi u)||_{L^p(Q)})
\crd
&\le 
A(||\n^mu||_{L^p(Q)}
+||\n^mw||_{L^p(Q)}).
\crd
&\le 
A||\n^mu||_{L^p(Q)}.
\cr}
$$

Here the triangle inequality has been used as well as weak Poincar\'e
inequalities.
\bs

Hence when the (*) {\bf Assumptions} are at hand it suffices for
(3.2) to hold for Spectral synthesis to hold as well.
\bs

The treatment of this problem has been taylored to suit the use of
polynomial capacities. However in this  situation 
only the higher order term in RHS of the Poincar\'e inequality (for a
cube) matters. Hence the $\G$ or $\T$ capacities are
are equivalent. 

However the previous calculations makes $\T$-capacities easier to handle.

We recall the definition

$$
\eq{
T^{\a}_{m,k,p,{\Cal A}}
=
\big\{
u\in {\Cal A}:
||\Pi_{m-1,k,p}u||_{L^p(Q)}
\ge
&\a||\Pi^{\text{\rm co}}_{m-1,k,p}u||_{L^p(Q)}
\crd
&\text{and}
\crd
||u-\Pi_{m-1,m-1,p}u||_{L^p(Q)}\le&{1\over 2}||\Pi_{m-1,m-1,p}u||_{L^p(Q)}
\big\}
\cr}
$$

\no
(i.e. in this case $\a$ is of no significance.) 

$$
\T^{\a}_{m,k,p}({\Cal A})
=\min \left[\inf_{u\in T^{\a}_{m,k,p,{\Cal A}}}
{||\n^mu||^p_{L^p(Q)}\over ||\Pi_{m-1,k,p} u||^p_{L^p(Q)}},1\right].
$$

Hence (3.2) now is formulated as a theorem. This theorem happens to be
very useful in a general sense as we hope to show in some later
paper(s). Here we make some consequences only.
\bs

\proclaim{3.3 Theorem} 
Spectral synthesis for Sobolev space with
parameters $m$, $p$ holds, i.e. for all closed $K$ in \R the two
subspaces $W^{m,p}_0(K^c)$ and $W^{m,p}_{0,m-1}(K^c)$ are equal 

iff for $Q_0$ a unit cube, for all $K$ closed, $K\subset\bar Q_0$
and $P\in{\Cal P}_{m-1}$

$$
\eq{
&\T_{m,k,p}(\{u\in W^{m,p}_0(K^c)|_{\bar Q_0}: \Pi u=P\})
\crd
&\le 
\crd
&A\cdot\T_{m,k,p}(\{u\in W^{m,p}_{0,m-1}(K^c)|_{\bar Q_0}: \Pi u=P\})
\cr}
\leqno(3.3)
$$
\endproclaim
\bs

\demo{Proof} First note that (3.3) is only a reformulation of (3.2)
since that formulation is dilation invariant when all $K$ are considered.

Before we do the remaining (*) parts of the proof we state a corollary
and a very analogous theorem.
\bs

\proclaim{3.4 Theorem} Spectral synthesis for Sobolev space with
parameters $1$, $p$ holds for $p\ge 1$.
\endproclaim
\bs

This is the earlier mentioned result by Hedberg.
\bs

\demo{Proof} This follows from Theorem 3.3 and Proposition 2.22 since
there really is an equivalence at hand in this Proposition.
\bs

The following is a new result. Its difficulty stems
from the fact that the Sobolev functions involved are not $C^1$
functions when $p\le N$, the dimension. The same argument applies of
course to the buissness of Sobolev space synthesis as a whole.
\bs

\proclaim{3.5 Theorem} For any closed $K$ in \R\ and $p\ge 1$

$$
W^{2,p}_{0,0}(K^c)_+=W^{2,p}_{0}(K^c)_+
$$

\no
holds.
\endproclaim
\bs

\demo{Proof} Apply the methods above then use Theorem 2.27 etc.
\bs

Now we turn to the points of the argument that have been postponed.
\bs

\demo{Proof {\rm of the (*) marked assumptions}}
\bs

As said before 
\bs

(i) $A_Q=A_0$ for all cubes $Q$
\bs

follows when the proper polynomial capacity equivalence is at hand as
is either proved or postulated in the treatment above.
\bs

(ii) The cubes cover $K$ at most a finite number $G$ times.
\bs

(iii) The property $O_j\subset O_{j+1}$.
\bs

These last two assumptions are taken care of by proceeding slightly
more delicately. 

Starting with the cubic lattice decomposition $\{Q\}$ then do the
redefining inside every $Q$. Next shift the lattice ${1\over 2}l(Q)$ in
some orthogonal direction.

However now we have to adjust to this new situation and define the
sets of functions within the cubes differently so that the
already accomplished redefinition stays when the new redefinition
takes place. Hence the good part grows as it should. 

This does not affect the polynomial capacity equivalences.

After iterating $N$ times we get that $O_r$ covers $K$. Hence $G=N$
and also (iii) is satisfied.
\bs

End of proof of the results for the whole section.


\Refs

\ref
\no[A-P]
\by D.R. Adams and J.C. Polking
\paper The equivalence of two definitions of capacity 
\pages 529-534
\jour Proc. Am. Math. Soc. 
\vol 37
\yr 1973
\endref

\ref
\no [A-H]
\by D.R. Adams and L.I. Hedberg
\book
Function spaces and potential theory
\publ Springer Verlag
\yr 1996
\endref

\ref
\no[RAD]
\by R.A. Adams
\book Sobolev Spaces
\publ Academic press
\yr 1975
\endref

\ref
\no[HED]
\by L.I. Hedberg
\paper Spectral synthesis in Sobolev spaces,
and uniqueness of solutions of the Dirichlet problem
\jour Acta Mathematica
\vol 147 
\yr 1981
\pages 237--264
\endref

\ref
\no[HED2]
\by L.I. Hedberg
\paper Two approximation problems in function spaces
\jour Arkiv f\"or matematik
\vol 16 
\yr 1978 
\pages 51--81
\endref

\ref
\no[MAZ]
\by V.G. Maz$'$ja
\book Sobolev Spaces
\yr 1985 
\publ Springer
\endref

\ref
\no[M-H]
\by V.G. Maz$'$ja and V.P. Havin
\paper Non-linear potential theory
\jour Uspehi Mat. Nauk
\vol 27:6 
\yr1972
\pages 67--138
\moreref
\paperinfo English translation
\jour Rus. math. surv. \vol 27, No 6
\yr 1972\pages 71--148
\endref

\ref
\no[MEY]
\by N.G. Meyers
\paper Continuity of Bessel potentials
\jour Israel J. of Math.
\vol 11 
\yr 1972 
\pages 271--283
\endref

\ref
\no[MEY2]
\by N.G. Meyers
\paper A theory of capacities for potentials of functions in
Lebesgue classes
\jour Math. Scand.
\vol 26
\year 1970
\pages 255--292
\endref

\ref
\no[NYS] 
\by K. Nystr\"om 
\paper Thesis
\paperinfo Univ. Ume\aa
\yr 1995
\endref


\ref
\no[WAN]
\by A. Wannebo
\paper Hardy inequalities
\vol 109
\yr 1990
\jour Proc. AMS
\pages 85-95
\endref
\bye